\documentclass[11pt]{article}
\usepackage[a4paper,margin=27mm,headheight=14pt]{geometry}
\usepackage[T1]{fontenc}
\usepackage{lmodern,amsmath,amssymb,amsthm}
\usepackage{microtype,enumitem,fancyhdr,needspace}
\setlist[enumerate]{itemsep=3pt,topsep=5pt}
\allowdisplaybreaks[1]
\usepackage[colorlinks=true,linkcolor=blue,citecolor=blue,urlcolor=blue]{hyperref}
\hypersetup{pdftitle={Alternating chains and lattice strict singularity},pdfauthor={Ma Yicen},pdfsubject={Candidate proof of LSS equals DSS over the real field; independent review pending}}
\newtheorem{theorem}{Theorem}[section]
\newtheorem{lemma}[theorem]{Lemma}
\newtheorem{proposition}[theorem]{Proposition}
\newtheorem{corollary}[theorem]{Corollary}
\theoremstyle{remark}
\newtheorem{remark}[theorem]{Remark}
\numberwithin{equation}{section}
\newcommand{\norm}[1]{\left\lVert#1\right\rVert}
\newcommand{\LSS}{\mathrm{LSS}}
\newcommand{\DSS}{\mathrm{DSS}}
\newcommand{\nucl}{\nu}

\title{Alternating chains and lattice strict singularity}
\author{Ma Yicen}
\date{Consolidated manuscript, version 2.1\\31 October 2026}
\begin{document}
\maketitle
\begin{abstract}
We present an alternating-chain argument intended to show that every lattice
strictly singular operator from a real Banach lattice to a real Banach space
is disjointly strictly singular. A hypothetical counterexample yields two
LSS operators on a reflexive sequence lattice whose difference is bounded
below. Quotienting by the range of that difference gives a common map
bounded below on a positive tail cone. Finite chains of positive disjoint
vectors, with alternating small differences, then lead to incompatible
length estimates. We record the resulting operator and sequence
consequences and extend the construction to common postcompositions with
DSS errors. Independently of the proposed equality, we prove simultaneous
nuclear restriction results for norm-$\sigma$-compact DSS families and
uniform nuclear tails on each compact piece. For norm-compact families
we also obtain common Hilbert-space factors and uniformly small nuclear
corrections having a common infinite-dimensional kernel. Examples identify
limits of these selection statements. The main equality is a candidate
real-field proof; no complex-field equality is asserted.
\end{abstract}


\medskip
\noindent\textbf{Manuscript status.}
This consolidated article contains a candidate complete proof over the real
field and its proposed consequences. Internal checks are not independent
review or formal verification. The general equality is presented for
mathematical assessment, not as an externally established resolution.
The DSS selection results in Sections~\ref{sec:dss}--\ref{sec:kernel}
have separate proofs valid over either scalar field. Their LSS formulations
invoke the proposed main theorem. Independent novelty and priority have not
been established.

\section{Introduction}

Let $E$ be a real Banach lattice and let $X$ be a real Banach space.
A bounded linear operator $T:E\to X$ is \emph{lattice strictly singular}
(LSS) if its restriction to every infinite-dimensional closed sublattice
fails to be bounded below. It is \emph{disjointly strictly singular}
(DSS) if it is bounded below on the closed span of no infinite sequence
of nonzero pairwise disjoint vectors. Every infinite-dimensional Banach
lattice contains an infinite positive disjoint sequence, so
\[
 \DSS(E,X)\subseteq\LSS(E,X).
\]
The reverse inclusion requires attention to signs: the span of a positive
disjoint sequence is a sublattice, whereas the span of an arbitrary
disjoint sequence need not be one.

Flores, L\'opez-Abad and Tradacete \cite{FLT} studied the converse and
related it to additivity of LSS operators and unconditional sequence
decompositions. Their results include equality for targets with the
unconditional subsequence property. Flores, Hern\'andez and Raynaud
\cite{FHR} had established equality of the corresponding super classes,
which impose uniform finite-dimensional conditions. The argument proposed
here addresses the ordinary classes over the real field.

\begin{theorem}[Proposed main theorem]\label{thm:main}
For every real Banach lattice $E$ and every real Banach space $X$,
\[
 \LSS(E,X)=\DSS(E,X).
\]
\end{theorem}

The proof has two components. A hypothetical counterexample first yields
two LSS operators $A,B:U\to X$, with $U$ a reflexive sequence lattice,
whose difference $J=A-B$ is bounded below. The range $J(U)$ is closed
and has an unconditional basis. Quotienting by this range supplies, on a
tail, a common map $C=QA=QB$ which is bounded below on positive vectors.
Only the particular subspace $J(U)$ is required to have an unconditional
basis; no such condition is imposed on $X$.

The second component is a finite construction. Small signed blocks for
$A$ and $B$ are assembled into alternating chains of positive disjoint
nodes whose $C$-images are asymptotically equal. A doubling procedure
gives chains of every dyadic length. A lower bound obtained from $C$
makes the norm of their alternating sum grow with the number of nodes.
A telescoping estimate for $J$ bounds that same norm using only the two
endpoints. The contradiction occurs at a fixed finite length.

Section~\ref{sec:prelim} records the reduction and selection facts.
Section~\ref{sec:quotient} proves the quotient cone estimate.
Section~\ref{sec:chains} contains the doubling construction, and
Section~\ref{sec:endpoint} gives the endpoint estimate and the proof of
the main theorem, together with a finite error threshold.
Section~\ref{sec:consequences} records operator and sequence consequences.
Sections~\ref{sec:factor}--\ref{sec:graph} give the common-factor extension.
Sections~\ref{sec:dss}--\ref{sec:kernel} establish independent DSS selection
and perturbation results, including a common nuclear restriction for
norm-$\sigma$-compact families and norm-continuous families parametrized
by $\mathbb R^d$. The final section records sharpness examples,
scope, and the main dependencies requiring review.

\section{Preliminaries and a reflexive reduction}\label{sec:prelim}

Throughout, all scalar fields are real, all operators are bounded and
linear, and $[v_n]$ denotes the closed linear span of a sequence.
An operator is \emph{bounded below} on a subspace when
$\norm{Tv}\geq c\norm v$ there for some $c>0$. For a vector $x$
in a Banach lattice, $x^+=x\vee0$, $x^-=(-x)\vee0$, and
$|x|=x^++x^-$. Disjointness means $|x|\wedge|y|=0$.

For positive disjoint vectors $v_j$ in a Banach lattice,
\begin{equation}\label{eq:modulus}
 \left|\sum_j a_jv_j\right|=\sum_j|a_j|v_j.
\end{equation}
In particular, if $\norm{v_j}=1$, then
\begin{equation}\label{eq:coeff}
 |a_j|\leq\norm{\sum_i a_iv_i}.
\end{equation}
The closed span of the $v_j$ is a closed sublattice.

We will work on a sequence Banach lattice $U$ whose coordinate vectors
$(e_n)$ form a normalized $1$-unconditional Schauder basis. A positive block
is a finite linear combination of the $e_n$ with nonnegative coefficients.
Write $U^{(m)}=[e_n:n\geq m]$. All tails have the inherited norm and
order. A successive block sequence has its finite supports ordered
from left to right; mere disjointness of supports does not require this
ordering. Normalized block sequences in
a reflexive such $U$ are weakly null.

For a normalized disjoint sequence, also over the complex field, its
coordinate functionals have norm one and its initial-segment projections
are contractive. The coefficient estimate remains valid for signed
disjoint vectors, since their moduli have disjoint supports.

An operator $S:F\to Y$ is nuclear if it has a representation
$Sx=\sum_n f_n(x)y_n$ with $\sum_n\norm{f_n}\norm{y_n}<\infty$.
Its nuclear norm $\nucl(S)$ is the infimum of these sums.
No approximation property is needed for the estimates used here.

\begin{lemma}[Summable images]\label{lem:nuclear}
Let $(v_n)$ be normalized and disjoint, and let $T:E\to Y$ be bounded.
If $\sum_n\norm{Tv_n}<\infty$, then $T|_F$ is nuclear for $F=[v_n]$, and
\begin{equation}\label{eq:nucleartail}
 \nucl\bigl(T|_F(I_F-P_m)\bigr)
 \leq\sum_{n>m}\norm{Tv_n}.
\end{equation}
Here $P_0=0$.
\end{lemma}
\begin{proof}
Use the representation
$T|_F=\sum_n v_n^*\otimes Tv_n$ and \eqref{eq:coeff}.
The series converges in nuclear norm, and agrees with $T$ on a dense
linear span, hence on $F$. Its tail gives \eqref{eq:nucleartail}.
\end{proof}

\begin{lemma}[Composition]\label{lem:composition}
If $T:E\to X$ is LSS and $R:V\to E$ is a bounded lattice homomorphism,
then $TR$ is LSS. A bounded postcomposition of an LSS operator is also LSS.
\end{lemma}
\begin{proof}
If $TR$ is bounded below on an infinite-dimensional closed sublattice $H$,
then $R|_H$ is bounded below, so $R(H)$ is an infinite-dimensional closed
sublattice. The inequalities
$\norm{TRh}\geq c\norm h\geq(c/\norm R)\norm{Rh}$ contradict LSS of $T$.
The zero-operator cases are immediate. The postcomposition assertion follows
directly from boundedness of the postcomposing operator.
\end{proof}

\Needspace{8\baselineskip}
\begin{lemma}[Small signed blocks]\label{lem:small}
Let $D:U\to X$ be LSS and let $(v_n)$ be positive, pairwise disjoint,
finitely supported vectors with $\inf_n\norm{v_n}\geq b>0$.
Given any sequence $\eta_k>0$, there are successive finite sets $I_k$ of
row indices and coefficients $a_n$ such that
\[
 z_k=\sum_{n\in I_k}a_nv_n,\qquad
 \norm{z_k}=1,\qquad \norm{Dz_k}<\eta_k,
 \qquad |a_n|\leq b^{-1}.
\]
The coefficients may have both signs.
\end{lemma}
\begin{proof}
Each tail spans an infinite-dimensional closed sublattice. The restriction
of $D$ is not bounded below there. Finite-support approximation, followed by
normalization, gives $z_k$ of the required size; repeat beyond the previous
row indices. The coefficient bound follows from
$|a_n|\norm{v_n}\leq\norm{z_k}$ and $\norm{z_k}=1$.
\end{proof}

\begin{lemma}[Reflexive two-column reduction]\label{lem:reduction}
If an LSS operator $T:E\to X$ is not DSS, there exist a reflexive sequence
lattice $U$ as above and LSS operators $A,B:U\to X$ such that
\begin{equation}\label{eq:J}
 J=A-B,\qquad \norm{Ju}\geq c\norm u\quad(u\in U)
\end{equation}
for some $c>0$.
\end{lemma}
\begin{proof}
Choose a normalized disjoint sequence $(x_n)$ on whose closed span $T$ is
bounded below, and complete $c_{00}$ in the norm
$\norm{a}_U=\norm{\sum_n a_nx_n}_E$. Disjointness makes this norm
$1$-unconditional. Set $R_\pm e_n=x_n^\pm$.
These maps are bounded lattice homomorphisms of norm at most one, since
\[
 \left|\sum_n a_nx_n^\pm\right|
 \leq\sum_n|a_n|\,|x_n|
 =\left|\sum_n a_nx_n\right|.
\]
Thus $A=TR_+$ and $B=TR_-$ are LSS by Lemma~\ref{lem:composition}, and
$Je_n=Tx_n$, giving \eqref{eq:J}.

We recall the established exclusion used in \cite[Proposition 2.1 and the
proof of Theorem 2.5]{FLT}: if an operator fixes a disjointly generated
copy of $c_0$ or $\ell_1$, it fixes an infinite-dimensional sublattice.
If $U$ were nonreflexive, the classical unconditional-basis criterion
would supply a block sequence equivalent to the unit vector basis of
$c_0$ or $\ell_1$. The corresponding blocks of $(x_n)$ are disjoint
and are still fixed by $T$, a contradiction. This is the reflexive
reduction in \cite[proof of Theorem 2.5]{FLT}; the underlying basis
criterion is due to James \cite{James}.
\end{proof}

\section{Positive-cone estimates in a quotient}\label{sec:quotient}

\begin{lemma}[Tail cone bound]\label{lem:cone}
Suppose $A,B$ are LSS and \eqref{eq:J} holds. On some coordinate tail,
$A$ is bounded below on the positive cone.
\end{lemma}
\begin{proof}
If this fails, positive finite-support approximation gives normalized
positive successive blocks $v_k$ such that
$\norm{Av_k}<\eta_k$ and $\sum_k\eta_k<c/2$. For a finite combination
$v=\sum_k a_kv_k$, \eqref{eq:coeff} gives
\[
 \norm{Av}\leq\Big(\sum_k\eta_k\Big)\norm v.
\]
Therefore $\norm{Bv}=\norm{Av-Jv}\geq(c/2)\norm v$.
By continuity this holds on the infinite-dimensional closed sublattice
$[v_k]$, contradicting LSS of $B$.
\end{proof}

A Banach space has the \emph{unconditional subsequence property}
(USP) if every seminormalized weakly null sequence has an unconditional
basic subsequence. We recall the selection argument used below.

\begin{lemma}[Unconditional subsequences]\label{lem:usp}
A Banach space with an unconditional Schauder basis has USP.
\end{lemma}
\begin{proof}
Let $(y_n)$ be a seminormalized weakly null sequence in such a space $Y$,
and put $a=\inf_n\norm{y_n}>0$. Write $P_m$ for the basis projections
and $K_0$ for an unconditional constant of the basis. Since
$P_my_n\to0$ for every fixed $m$, a gliding-hump selection gives a
subsequence $(y_{n_k})$ and successive basis blocks $(w_k)$ with
$\sum_k\norm{y_{n_k}-w_k}$ arbitrarily small: at each step choose the
next index to make the previous initial segment small, then choose the
new support endpoint to make the remaining tail small.

Arrange $\norm{w_k}\geq a/2$ and
$\sum_k\norm{y_{n_k}-w_k}<a/(4K_0)$. Coordinate suppression implies
\[
 |t_k|\leq\frac{2K_0}{a}\norm{\sum_i t_iw_i}.
\]
Thus, on $[w_k]$, the map $\sum t_kw_k\mapsto\sum t_ky_{n_k}$
differs from the inclusion into $Y$ by an operator of norm less than
$1/2$. It is therefore an isomorphism onto its closed range.
Since $(w_k)$ is unconditional, so is $(y_{n_k})$.
This is the classical block-selection argument; see also \cite{BP}.
\end{proof}

\begin{lemma}[Passage to a quotient]\label{lem:quotient}
Let $U$ be reflexive, and let $D:U\to X$ be LSS with
\[
 \norm{Du}\geq d\norm u\quad(u\in U_+),\qquad d>0.
\]
Let $Y\subseteq X$ be closed with USP, and let $Q:X\to X/Y$ be the
quotient map. Then $QD$ is bounded below on the positive cone of a tail.
\end{lemma}
\begin{proof}
Otherwise, using positive finite-support approximation, there are
normalized positive successive blocks $v_n$ with
$\norm{QDv_n}\to0$. Choose $y_n\in Y$ so that
$e'_n=Dv_n-y_n\to0$ in norm. Since $U$ is reflexive, $(v_n)$ is weakly
null. Hence $(y_n)$ is weakly null in $X$ and also in $Y$, by
Hahn--Banach. After discarding finitely many terms it is seminormalized,
since $\norm{Dv_n}\geq d$ and $\norm{Dv_n}\leq\norm D$.

Pass to an unconditional subsequence of $(y_n)$, of unconditional constant
$K$, and then thin further to arrange
$\eta=\sum_n\norm{e'_n}<d/[2(K+1)]$.
For $v=\sum_n a_nv_n$ and $w=\sum_n|a_n|v_n$, disjointness gives
$\norm v=\norm w$ and $\sup_n|a_n|\leq\norm v$. Therefore
\begin{align*}
 \norm{Dv}
 &\geq \norm{\sum_n a_ny_n}-\eta\norm v\\
 &\geq K^{-1}\norm{\sum_n|a_n|y_n}-\eta\norm v\\
 &\geq \left(\frac dK-\left(1+\frac1K\right)\eta\right)\norm v
 >\frac d{2K}\norm v.
\end{align*}
Thus $D$ fixes the infinite-dimensional closed sublattice $[v_n]$, a
contradiction.
\end{proof}

\begin{proposition}[A common quotient map]\label{prop:common}
For the operators in Lemma~\ref{lem:reduction}, one may pass to a
coordinate tail and find a Banach space $Z$, a contraction $Q:X\to Z$
and a constant $\delta>0$ such that
\begin{equation}\label{eq:common}
 C=QA=QB:U\longrightarrow Z,\qquad
 \norm{Cu}\geq\delta\norm u\quad(u\in U_+).
\end{equation}
The lower bound $c$ for $J$ remains valid on the tail.
\end{proposition}
\begin{proof}
Write $U_0$ for the domain before taking tails and set $Y=J(U_0)$.
Because $J$ is bounded below, $Y$ is closed and isomorphic to $U_0$.
It therefore has an unconditional basis and USP by Lemma~\ref{lem:usp}.
Put $Z=X/Y$ and let $Q$ be the canonical quotient map.

Lemma~\ref{lem:cone} supplies an $A$-cone estimate on a tail.
Apply Lemma~\ref{lem:quotient} to the restriction of $A$ to that tail.
On a further tail it gives the required lower bound for $QA$.
Since $QJ=0$, we have $QA=QB$. Relabel the final tail as $U$.
All restrictions remain LSS by Lemma~\ref{lem:composition}.
\end{proof}

\begin{remark}\label{rem:quotient}
The quotient is by the fixed subspace $J(U_0)$, even after taking
a smaller domain. Neither $X$ nor $Z$ is assumed to have USP.
No order structure or positivity is imposed on $A,B,C$, or $Q$.
The adjective \emph{positive} refers only to vectors in the domain.
\end{remark}

From now on, $c$ and $\delta$ are fixed. Write $L=\norm C>0$ and
$M=\norm B$.

\section{Alternating chains of disjoint positive vectors}\label{sec:chains}

The construction in this section uses the common map in
\eqref{eq:common} but not reflexivity. It produces a different infinite
array for each prescribed finite width.

\begin{lemma}[Arbitrarily long finite chains]\label{lem:chains}
Assume $A,B:U\to X$ are LSS and \eqref{eq:common} holds for a quotient
map $Q$ of norm at most one. For every $N=2^k$, $k\geq1$, there is an
array of positive finitely supported vectors
\[
 (u_{n,j})_{n\geq1,\ 1\leq j\leq N}
\]
which are pairwise disjoint over all pairs $(n,j)$, with
$\norm{u_{n,1}}=1$, such that, as $n\to\infty$,
\begin{align}
 \max_{1\leq j\leq N}\norm{Cu_{n,j}-Cu_{n,1}}&\longrightarrow0,
                                                        \label{eq:Cchain}\\
 \max_{1\leq i\leq N/2}\norm{A(u_{n,2i-1}-u_{n,2i})}&\longrightarrow0,
                                                        \label{eq:Achain}\\
 \max_{1\leq i<N/2}\norm{B(u_{n,2i}-u_{n,2i+1})}&\longrightarrow0.
                                                        \label{eq:Bchain}
\end{align}
The maximum over the empty set in \eqref{eq:Bchain}, for $N=2$, is zero.
\end{lemma}
\begin{proof}
Let $\varepsilon_n$ denote the maximum of the three errors for a row.
Whenever an array has been constructed, \eqref{eq:common} and
$\norm{u_{n,1}}=1$ imply
\begin{equation}\label{eq:uniformnodes}
 \frac{\delta-\varepsilon_n}{L}\leq\norm{u_{n,j}}
 \leq\frac{L+\varepsilon_n}{\delta}.
\end{equation}
Thus all columns are uniformly bounded above and away from zero after
discarding finitely many rows.

\textit{Base step.} Apply Lemma~\ref{lem:small} to $A$ and $(e_n)$ to
obtain disjoint normalized signed blocks $z_n$ with $\norm{Az_n}\to0$.
Write $p_n=z_n^+$ and $q_n=z_n^-$. Then
$\norm{p_n+q_n}=1$ and
$\norm{Cp_n-Cq_n}=\norm{QAz_n}\to0$.
Since $p_n+q_n\geq0$,
\[
 \delta\leq\norm{Cp_n+Cq_n}
 \leq 2\norm{Cp_n}+\norm{Cp_n-Cq_n}.
\]
In particular $\norm{p_n}\geq\delta/(4L)$ eventually. Set
$u_{n,1}=p_n/\norm{p_n}$ and $u_{n,2}=q_n/\norm{p_n}$.
All required errors tend to zero, and disjointness is preserved.

\textit{Doubling step.} Suppose the assertion holds for an even $N$.
Pass to a subsequence of rows, keeping the notation, such that
$\varepsilon_n\leq2^{-n}$ and $\varepsilon_n\leq\delta/2$.
By \eqref{eq:uniformnodes}, the last column satisfies
$\norm{u_{n,N}}\geq b:=\delta/(2L)$.
Apply Lemma~\ref{lem:small} to $B$ on this positive disjoint column:
\[
 z_k=\sum_{n\in I_k}a_nu_{n,N},\qquad
 \norm{z_k}=1,\qquad \beta_k:=\norm{Bz_k}\longrightarrow0,
 \qquad |a_n|\leq b^{-1}.
\]
For each old column $j$, use exactly the same coefficients to define
\[
 P_{k,j}=\sum_{\substack{n\in I_k\\a_n>0}}a_nu_{n,j},\qquad
 R_{k,j}=\sum_{\substack{n\in I_k\\a_n<0}}(-a_n)u_{n,j}.
\]
All these vectors are positive and pairwise disjoint, including between
different $k$. Indeed, different columns use disjoint supports, and the
positive and negative index sets within each column are disjoint. Put
\[
 t_k=b^{-1}\sum_{n\geq\min I_k}\varepsilon_n\longrightarrow0.
\]
For example, each inherited $A$-edge in the $P$-path satisfies
\[
 \norm{A(P_{k,2i-1}-P_{k,2i})}
 \leq\sum_{\substack{n\in I_k\\a_n>0}}|a_n|\varepsilon_n
 \leq t_k.
\]
The same estimate controls every inherited $A$- or $B$-edge in both
paths, independently of the number of indices in $I_k$. Moreover
\begin{align*}
 \norm{CP_{k,j}-CP_{k,1}}&\leq t_k,&
 \norm{CR_{k,j}-CR_{k,1}}&\leq t_k,\\
 \norm{B(P_{k,N}-R_{k,N})}&=\beta_k,&
 \norm{CP_{k,N}-CR_{k,N}}&\leq\beta_k.
\end{align*}
The last inequality uses $C=QB$. All $C$-images of the $2N$ new vectors
are consequently within $3t_k+\beta_k$ of $CP_{k,1}$.

Join the two paths at their last endpoints and reverse the second path:
\begin{equation}\label{eq:join}
 P_{k,1},\ldots,P_{k,N},R_{k,N},R_{k,N-1},\ldots,R_{k,1}.
\end{equation}
The joining edge is a $B$-edge. Since $N$ is even, the full list again
has $A$-edges in positions $(1,2),(3,4),\ldots$, and $B$-edges in
positions $(2,3),(4,5),\ldots$. Its first and last edges are $A$-edges.

It remains to justify a uniform normalization. Here
$p=P_{k,N}$ and $r=R_{k,N}$ are precisely the positive and negative
parts of $z_k$. Thus $\norm{p+r}=1$, and
\[
 \delta\leq\norm{Cp+Cr}\leq2\norm{Cp}+\beta_k,
 \qquad \norm{Cp}\geq(\delta-\beta_k)/2.
\]
Also $\norm{CP_{k,1}-Cp}\leq t_k$. Therefore, after discarding finitely
many $k$,
\begin{equation}\label{eq:normalization}
 \norm{P_{k,1}}\geq\frac{\norm{CP_{k,1}}}{L}
 \geq\frac{\delta-\beta_k-2t_k}{2L}
 \geq\frac{\delta}{4L}.
\end{equation}
Divide the entire list \eqref{eq:join} by $\norm{P_{k,1}}$.
The first vector now has norm one, and the maximum of the three new
errors is bounded by
\begin{equation}\label{eq:newerror}
 \frac{4L}{\delta}(3t_k+\beta_k)\longrightarrow0.
\end{equation}
Disjointness is unchanged. This constructs the array of width $2N$
and completes the finite induction.
\end{proof}

\begin{remark}
The columns inside a row need not be successive in coordinate order.
They are pairwise disjoint, which is exactly what the proof uses.
The near-kernel extraction is always applied to an infinite positive
disjoint column. Signed blocks are split into their positive and negative
parts before the next application of LSS. For each prescribed finite $N$,
only finitely many doubling steps are made.
\end{remark}

\section{The endpoint bound and the main theorem}\label{sec:endpoint}

The obstruction to long chains is a finite estimate. In particular,
it does not involve an infinite chain or the extraction of an infinite
$\ell_1$ sequence from finite-dimensional estimates.

\begin{theorem}[Chain-length bound]\label{thm:length}
Suppose \eqref{eq:J} and \eqref{eq:common} hold. Let
$(u_{n,j})_{n\geq1,\,1\leq j\leq N}$ be an array of positive vectors,
pairwise disjoint within each row, where $N\geq2$ is even and
$\norm{u_{n,1}}=1$. If \eqref{eq:Cchain}--\eqref{eq:Bchain} hold, then
\begin{equation}\label{eq:length}
 N\leq\frac{2\norm B\,\norm C}{c\delta}.
\end{equation}
\end{theorem}
\begin{proof}
Fix a row, suppress its index, and write $\varepsilon$ for its maximum
error and $\tau=\norm{Cu_1}$. Then $\delta\leq\tau\leq L$.
Set
\[
 z=\sum_{j=1}^{N}(-1)^ju_j,
 \qquad w=\sum_{j=1}^{N}u_j=|z|.
\]
Using $C$ to estimate the norm of the positive sum gives
\begin{equation}\label{eq:long}
 \norm z=\norm w
 \geq\frac{\norm{Cw}}{L}
 \geq\frac{N\tau-(N-1)\varepsilon}{L}.
\end{equation}
The alternating edge estimates give
\begin{align*}
 Az&=\sum_{i=1}^{N/2} A(u_{2i}-u_{2i-1}),\\
 Bz&=Bu_N-Bu_1+
       \sum_{i=1}^{N/2-1}B(u_{2i}-u_{2i+1}).
\end{align*}
Consequently, since $J=A-B$,
\begin{equation}\label{eq:short}
 c\norm z\leq\norm{Jz}
 \leq M(\norm{u_1}+\norm{u_N})+(N-1)\varepsilon
 \leq\frac{M(2\tau+\varepsilon)}{\delta}+(N-1)\varepsilon.
\end{equation}
The last inequality uses $\norm{u_1}=1\leq\tau/\delta$ and
$\norm{u_N}\leq\norm{Cu_N}/\delta\leq(\tau+\varepsilon)/\delta$.
Combining \eqref{eq:long} and \eqref{eq:short}, and dividing by
$\tau\geq\delta$, yields
\begin{equation}\label{eq:finitebound}
 \frac{cN}{L}
 \leq\frac{2M}{\delta}
   +\left(\frac{c(N-1)}{L}+\frac M\delta+N-1\right)
     \frac{\varepsilon}{\tau}.
\end{equation}
Let the row index tend to infinity. The error term tends to zero, proving
\eqref{eq:length}.
\end{proof}

\begin{corollary}[Finite error threshold]\label{cor:threshold}
Let $U$ be a real Banach lattice, let $A,B:U\to X$ and $C:U\to Z$
be bounded, and suppose that $A-B$ is bounded below by $c>0$ and
$\norm{Cu}\geq\delta\norm u$ for $u\geq0$. Set $L=\norm C$ and
$M=\norm B$. For positive pairwise disjoint $u_1,\ldots,u_N$, where
$N\geq2$ is even and $\norm{u_1}=1$, put $\tau=\norm{Cu_1}$ and
\[
 \varepsilon=\max\left\{
 \max_j\norm{Cu_j-Cu_1},\;
 \max_{1\leq i\leq N/2}\norm{A(u_{2i-1}-u_{2i})},\;
 \max_{1\leq i<N/2}\norm{B(u_{2i}-u_{2i+1})}\right\}.
\]
If $N>2ML/(c\delta)$, then
\begin{equation}\label{eq:threshold}
 \frac{\varepsilon}{\tau}\geq
 \frac{c\delta N-2ML}{(c+L)\delta(N-1)+ML}>0.
\end{equation}
\end{corollary}
\begin{proof}
The single-row calculation \eqref{eq:long}--\eqref{eq:finitebound}
uses only the hypotheses just stated. Multiply \eqref{eq:finitebound}
by $L\delta$ and rearrange. The denominator is positive.
\end{proof}

\begin{remark}
Corollary~\ref{cor:threshold} does not assume LSS, reflexivity, or the
factorization $C=QA=QB$. It is a finite geometric estimate independent
of the proposed main theorem. For fixed operators its right-hand side
tends to $c/(c+L)$ as even $N$ tends to infinity. This does not give a
uniform quantitative comparison of LSS and DSS: the quotient argument
provides no uniform control of $\delta$ in terms of the original operator.
\end{remark}

\begin{proof}[Proof of Theorem~\ref{thm:main}]
Suppose $T:E\to X$ is LSS but not DSS. Lemma~\ref{lem:reduction} gives
$U,A,B,J,c$ with $U$ reflexive. Proposition~\ref{prop:common}
gives a tail on which \eqref{eq:common} holds. Restrict all operators
to this tail and fix their constants.

Choose a power of two $N>2\norm B\norm C/(c\delta)$.
Lemma~\ref{lem:chains} constructs an array of width $N$, whereas
Theorem~\ref{thm:length} forbids that width. This contradiction proves
$\LSS\subseteq\DSS$. The reverse inclusion was noted above.
\end{proof}

\section{Operator and sequence consequences}\label{sec:consequences}

The results in this section which invoke Theorem~\ref{thm:main} inherit
its candidate status. They are not independent evidence for that theorem.

\begin{corollary}[Additivity and norm closure]\label{cor:additivity}
For every real Banach lattice $E$ and real Banach space $X$,
$\LSS(E,X)$ is a norm-closed linear subspace of $\mathcal L(E,X)$.
In particular, an operator-norm-convergent series of LSS operators
has an LSS sum.
\end{corollary}
\begin{proof}
The DSS class is additive, as recalled in \cite[Introduction]{FLT}.
Theorem~\ref{thm:main} therefore gives additivity of LSS; scalar
multiplication is immediate. Norm closure follows directly from the
definition: if $T_n\to T$ and $T$ is bounded below by $a>0$ on an
infinite-dimensional closed sublattice $H$, then $T_n$ is bounded below
by $a/2$ on $H$ whenever $\norm{T_n-T}<a/2$.
Apply closure to the finite partial sums for the last assertion.
\end{proof}

\begin{remark}
Additivity is a known equivalent formulation of the LSS--DSS problem
\cite[Theorem 2.5]{FLT}, not a separate resolution of another problem.
Norm convergence cannot be replaced by strong operator convergence:
the finite-rank coordinate projections on $\ell_2$ are LSS and converge
strongly to the identity.
\end{remark}

\begin{corollary}[Lattices with a disjoint Schauder basis]\label{cor:basis}
If a real Banach lattice $E$ has a Schauder basis of pairwise disjoint
vectors, then, for every real Banach space $X$,
\[
 \LSS(E,X)=\DSS(E,X)=\mathrm{SS}(E,X),
\]
where $\mathrm{SS}$ denotes the strictly singular operators.
\end{corollary}
\begin{proof}
We give the usual block argument for $\DSS=\mathrm{SS}$; only the first
equality uses Theorem~\ref{thm:main}. If $T$ is not strictly singular,
choose an infinite-dimensional closed subspace $H$ on which $T$ is
bounded below by $a>0$. Write $P_m$ for the basis projections.
Recursively choose normalized $h_k\in H\cap\ker P_{m_{k-1}}$ and
$m_k>m_{k-1}$ so that $v_k=P_{m_k}h_k$ satisfies
$\norm{h_k-v_k}<\eta_k$. The vectors $v_k$ are disjoint blocks.
Choose $\eta_k<1/2$ and $s=2\sum_k\eta_k$ so small that
$a-(a+\norm T)s>0$. Disjointness and $\norm{v_k}\geq1/2$ give
$|t_k|\leq2\norm{\sum_i t_iv_i}$. For $v=\sum t_kv_k$ and
$h=\sum t_kh_k$, we have $\norm{h-v}\leq s\norm v$, and hence
\[
 \norm{Tv}\geq a\norm h-\norm T\norm{v-h}
 \geq\bigl(a-(a+\norm T)s\bigr)\norm v.
\]
Thus $T$ is not DSS. The opposite implication is immediate.
\end{proof}

\begin{lemma}[Positive block approximation]\label{lem:positiveblocks}
Let $U$ be a sequence lattice with a $1$-unconditional Schauder basis,
and let $D:U\to X$ be bounded. If $D$ is not LSS, there are normalized
positive successive finite blocks $(v_k)$ on whose closed span $D$ is
bounded below.
\end{lemma}
\begin{proof}
Choose normalized positive pairwise disjoint $h_n$ in a sublattice on
which $D$ is bounded below by $a>0$. For each fixed $m$, at most $m$
of the $h_n$ have a nonzero coordinate among $1,\ldots,m$. Recursively
select $h_{n_k}$ vanishing on the preceding initial segment and truncate
its tail to obtain positive successive finite blocks $w_k$, with
$\norm{w_k-h_{n_k}}<\eta_k$. Take $\eta_k<1/2$ and
$s=2\sum\eta_k$ with $a-(a+\norm D)s>0$.
The estimate in the proof of Corollary~\ref{cor:basis}, with $D$ in
place of $T$, proves that $D$ is bounded below on $[w_k]$.
Normalize $w_k$ to obtain $v_k$ without changing the span.
\end{proof}

\begin{theorem}[Finite-component positive block selection]\label{thm:components}
Let $(z_n)$ be an unconditional basic sequence in a real Banach space
$X$, and let $m\geq1$. Suppose
\[
 z_n=\sum_{j=1}^m x_n^{(j)},\qquad
 \norm{\sum_n a_nx_n^{(j)}}\leq K_j\norm{\sum_n a_nz_n}
 \quad(1\leq j\leq m)
\]
for every finitely supported real sequence $(a_n)$.
Then there are a fixed $j_0$, successive finite sets $I_k$, and
coefficients $\alpha_n\geq0$ such that the nonzero blocks
\[
 b_k=\sum_{n\in I_k}\alpha_nz_n,\qquad
 d_k=\sum_{n\in I_k}\alpha_nx_n^{(j_0)}
\]
are equivalent basic sequences. In particular, $(d_k)$ is unconditional.
No basic-sequence assumption is required on any individual column
$(x_n^{(j)})_n$.
\end{theorem}
\begin{proof}
Give $c_{00}$ the lattice norm
\[
 \norm{a}_U=\sup_{|\theta_n|\leq1}
                 \norm{\sum_n\theta_na_nz_n}
\]
and take its completion. If $K$ is an unconditional constant of $(z_n)$,
then $\norm{Ja}\leq\norm a_U\leq K\norm{Ja}$ for
$Je_n=z_n$. Thus $J:U\to X$ is bounded below. The coordinate basis
of $U$ is $1$-unconditional; it need not be normalized.
The domination hypotheses define bounded maps $A_je_n=x_n^{(j)}$,
and $\sum_{j=1}^m A_j=J$.

If every $A_j$ were LSS, Corollary~\ref{cor:additivity} would make
$J$ LSS, a contradiction. Fix $j_0$ for which $A_{j_0}$ is not LSS.
Lemma~\ref{lem:positiveblocks} supplies positive successive finite
blocks $v_k=\sum_{n\in I_k}\alpha_ne_n$ on whose span $A_{j_0}$
is bounded below. Both $J$ and $A_{j_0}$ are isomorphisms on this
span onto their respective images. Their image sequences
$b_k=Jv_k$ and $d_k=A_{j_0}v_k$ are therefore equivalent, as claimed.
\end{proof}

\begin{remark}
The conclusion preserves the same coefficients in the chosen component
and the original sequence, with one component fixed for every block.
The equivalence constants may depend on the decomposition and the chosen
blocks; no uniform bound is asserted. This is a consequence of the main
theorem and its additivity formulation. Its novelty relative to existing
sequence-selection formulations has not been established.
\end{remark}

\subsection*{The two-column formulation}

For a sequence $(s_j)$, a positive block sequence means nonzero vectors
$v_k=\sum_{j\in F_k}t_js_j$, where $t_j\geq0$ and the finite sets
$F_k$ are successive. Write $\asymp$ for two inequalities whose constants
are independent of all finitely supported scalar sequences.

\begin{corollary}\label{cor:sequence}
Let $X$ be a real Banach space and let $(x_n)$ and $(y_n)$ be basic
sequences in $X$. Suppose $(x_n-y_n)$ is an unconditional basic sequence
and
\begin{equation}\label{eq:sequence}
 \norm{\sum_n a_n(x_n-y_n)}
 \asymp\max\left\{\norm{\sum_n a_nx_n},
                       \norm{\sum_n a_ny_n}\right\}
\end{equation}
for every finitely supported real scalar sequence $(a_n)$.
Then $(x_1,y_1,x_2,y_2,\ldots)$ admits an unconditional positive
block sequence.
\end{corollary}
\begin{proof}
If no such block sequence existed, conditions (b.1)--(b.3) of
\cite[Theorem 2.5]{FLT} would hold. That equivalence gives a real Banach
lattice $E$ and an LSS non-DSS operator $E\to X$, contradicting
Theorem~\ref{thm:main}.
\end{proof}

In the terminology of \cite[Definition 2.6]{FLT}, the same equivalence
gives the unconditional decomposition property for every real Banach
space. These sequence conclusions depend on the proposed main theorem;
they are not independent verifications of it.

\section{Two common postcompositions}\label{sec:factor}

\begin{theorem}[Common-factor assertion]\label{thm:factor}
Let $A:E\to X$ and $B:E\to Y$ be LSS. Suppose bounded operators
$R:X\to Z$ and $S:Y\to Z$ satisfy
\begin{equation}\label{eq:factor}
 C=RA=SB,\qquad \norm{Cu}\geq\delta\norm u\quad(u\in E_+),
 \qquad\delta>0.
\end{equation}
Then $(A,B):E\to X\oplus_\infty Y$ is LSS. Neither reflexivity of $E$
nor equality of $X$ and $Y$ is required.
\end{theorem}
\begin{proof}
Fix an infinite-dimensional closed sublattice $H$ and choose a normalized
positive disjoint sequence in $H$. Its closed span is a sequence lattice
$U$ with a normalized $1$-unconditional basis. Restrictions of $A$ and $B$
are LSS and \eqref{eq:factor} is preserved. Put $L=\norm{C|_U}>0$ and
$M=\norm{B|_U}$.

We describe explicitly how the chain construction of Lemma~\ref{lem:chains}
changes. For each dyadic even $N$, it produces positive finite blocks
$u_{n,j}\in U$, disjoint over all $(n,j)$, with $\norm{u_{n,1}}=1$ and
row errors tending to zero:
\begin{align}
 \max_j\norm{Cu_{n,j}-Cu_{n,1}}&\longrightarrow0,\label{eq:cchain}\\
 \max_i\norm{A(u_{n,2i-1}-u_{n,2i})}&\longrightarrow0,\label{eq:achain}\\
 \max_{i<N/2}\norm{B(u_{n,2i}-u_{n,2i+1})}&\longrightarrow0.\label{eq:bchain}
\end{align}
To start, choose normalized signed blocks $z_n$ with $\norm{Az_n}\to0$
using LSS on positive tails. For $p_n=z_n^+$ and $q_n=z_n^-$,
\[
 \norm{Cp_n-Cq_n}\leq\norm R\,\norm{Az_n}\longrightarrow0,\qquad
 \delta\leq\norm{Cp_n+Cq_n}.
\]
Thus $\norm{p_n}\geq\delta/(4L)$ eventually; normalize both parts by
$\norm{p_n}$.

For completeness, the doubling estimates are as follows. Thin an array of
width $N$ so its maximum errors satisfy
$\varepsilon_n\leq\min(2^{-n},\delta/2)$. Its last column has norms at least
$b=\delta/(2L)$. LSS of $B$ supplies successive finite row blocks
\[
 z_k=\sum_{n\in I_k}a_nu_{n,N},\quad \norm{z_k}=1,\quad
 \beta_k=\norm{Bz_k}\to0,\quad |a_n|\leq b^{-1}.
\]
Use the same coefficients in every column, writing
\[
 P_{k,j}=\sum_{\substack{n\in I_k\\a_n>0}}a_nu_{n,j},\qquad
 R_{k,j}=\sum_{\substack{n\in I_k\\a_n<0}}(-a_n)u_{n,j}.
\]
All these vectors are positive and disjoint. With
$t_k=b^{-1}\sum_{n\geq\min I_k}\varepsilon_n\to0$, each inherited edge error
and each within-path $C$-difference is at most $t_k$. The new bridge satisfies
\[
 \norm{B(P_{k,N}-R_{k,N})}=\beta_k,\qquad
 \norm{CP_{k,N}-CR_{k,N}}\leq\norm S\,\beta_k.
\]
The ordered list
$P_{k,1},\ldots,P_{k,N},R_{k,N},\ldots,R_{k,1}$
has the required alternating edge types. Every $C$-image lies within
$3t_k+\norm S\,\beta_k$ of $CP_{k,1}$. Moreover,
\[
 \norm{P_{k,1}}\geq
 \frac{\delta-\norm S\,\beta_k-2t_k}{2L}\geq\frac{\delta}{4L}
\]
eventually. Dividing the list by $\norm{P_{k,1}}$ makes all new errors
at most $(4L/\delta)(3t_k+(1+\norm S)\beta_k)\to0$.
This proves \eqref{eq:cchain}--\eqref{eq:bchain} by finite induction.

For a fixed width let
$z_n=\sum_{j=1}^N(-1)^ju_{n,j}$ and $\tau_n=\norm{Cu_{n,1}}\geq\delta$.
Disjointness and the common-image estimates give
\[
 \norm{z_n}\geq\frac{N\tau_n-(N-1)\varepsilon_n}{L}.
\]
The $A$-edges yield $\norm{Az_n}\leq N\varepsilon_n/2$, while telescoping
the $B$-edges gives
\[
 \norm{Bz_n}\leq
 \frac{M(2\tau_n+\varepsilon_n)}{\delta}
 +(N/2-1)\varepsilon_n.
\]
Consequently
\[
 \limsup_n\frac{\max(\norm{Az_n},\norm{Bz_n})}{\norm{z_n}}
 \leq\frac{2ML}{N\delta}.
\]
Choose $N$ first and then the row index. There are unit vectors in $H$
on which both operators are arbitrarily small. As $H$ was arbitrary,
the joint operator is LSS.
\end{proof}

\begin{corollary}[DSS errors]\label{cor:errors}
The conclusion of Theorem~\ref{thm:factor} remains valid if its factorization
is replaced by
\[
 C=RA+D_A=SB+D_B,
\]
where $D_A,D_B:E\to Z$ are DSS and $C$ has the same positive-cone lower bound.
\end{corollary}
\begin{proof}
The established perturbation lemma says that LSS plus DSS is LSS
\cite[Lemma 2.4]{FLT}. Thus $\widehat A=(A,D_A)$ and
$\widehat B=(B,D_B)$ are LSS: write each as the sum of the coordinate
injections of its components. Bounded postcomposition preserves both
classes. Set $\widehat R(x,z)=Rx+z$ and $\widehat S(y,z)=Sy+z$.
Theorem~\ref{thm:factor} applies, and projection onto the original
coordinates proves the assertion.
\end{proof}

\section{Connected families and nuclear restrictions}\label{sec:graph}

\begin{corollary}\label{cor:graph}
Let $A_i:E\to X_i$, $i\in I$, be LSS, with $I$ finite or countably infinite.
Suppose a connected graph on $I$ has, for each edge $\{i,j\}$, bounded
postcompositions and DSS errors such that
\[
 C_{ij}=R_{ij}A_i+D_{ij}
       =S_{ij}A_j+\widetilde D_{ij},\qquad
 \norm{C_{ij}u}\geq\delta_{ij}\norm u\quad(u\geq0),\quad\delta_{ij}>0.
\]
Then every finite subfamily has an LSS joint operator. In every
infinite-dimensional closed sublattice $H\subset E$, there is a normalized
disjoint sequence $(v_n)$ such that every $A_i|_{[v_n]}$ is nuclear.
\end{corollary}
\begin{proof}
For a finite connected set, add vertices along a spanning tree. Suppose the
joint map $K$ for the existing vertices is LSS. An edge from vertex $i$ in
that set to the next vertex $j$ gives the same factorization for $K$ and
$A_j$, using $R_{ij}$ after the $i$th coordinate projection.
Corollary~\ref{cor:errors} proves that the enlarged joint map is LSS.
Any finite set lies in a finite connected set, so projection gives the
first assertion.

Choose increasing finite connected sets $F_n$ whose union is $I$, eventually
constant if $I$ is finite. Fix a positive normalized disjoint sequence
$(e_k)$ in $H$. On each unused tail of $[e_k]$, the joint map for $F_n$
is not bounded below. Choose a normalized signed finite block $v_n$
beyond the preceding block so that
\[
 \max_{i\in F_n}\norm{A_iv_n}<2^{-n}.
\]
For fixed $i$, this estimate holds for all sufficiently large $n$.
Hence $\sum_n\norm{A_iv_n}<\infty$, and Lemma~\ref{lem:nuclear} applies.
\end{proof}

\begin{remark}
No uniform bound on the edge constants is needed. The sequence $(v_n)$
consists of signed blocks of a fixed positive disjoint sequence in $H$.
Its span need not be a sublattice. The conclusion is nuclearity of each
restriction, not a uniform nuclear-norm bound over an arbitrary countable
family.
\end{remark}

\section{DSS families without common-factor assumptions}\label{sec:dss}

This section and the subsequent DSS results work over either scalar field.
They do not use the chain construction or Theorem~\ref{thm:main}. The underlying
small-restriction method is classical; compare the compact-restriction
characterization recalled in \cite[Introduction]{FHKT}.

\begin{lemma}[Finite simultaneous selection]\label{lem:finite}
Let $T_j:E\to Y_j$, $1\leq j\leq q$, be DSS and let $(x_n)$ be normalized and
disjoint. In every tail of this sequence and for every $\epsilon>0$, there
is a normalized finite block $v$ with $\max_j\norm{T_jv}<\epsilon$.
\end{lemma}
\begin{proof}
For $T_1$, its failure to be bounded below on each disjoint tail permits
successive normalized finite blocks $w_k$ with
$\sum_k\norm{T_1w_k}<\epsilon$. Equation~\eqref{eq:coeff} gives
$\norm{T_1|_{[w_k]}}<\epsilon$. Repeat within this block sequence for $T_2$,
and then for each remaining operator. Every new sequence is still disjoint
and the preceding norm bounds persist on smaller spans. After finitely
many steps, take one block.
\end{proof}

\begin{theorem}[Countable simultaneous nuclear restrictions]\label{thm:countable}
Let $T_j:E\to Y_j$, $j\geq1$, be DSS. Given any normalized disjoint
sequence $(x_n)$, there is a normalized successive block sequence $(v_n)$
of $(x_n)$ such that every $T_j|_{[v_n]}$ is nuclear.
\end{theorem}
\begin{proof}
At stage $n$, use Lemma~\ref{lem:finite} on the unused tail to arrange
$\norm{T_jv_n}<2^{-n}$ for $j\leq n$. For fixed $j$, the remaining
$j-1$ image norms form a finite sum. Apply Lemma~\ref{lem:nuclear}.
\end{proof}

\begin{corollary}[Norm-separable families]\label{cor:separable}
If $\mathcal A\subset\DSS(E,Y)$ is separable in operator norm, there is
a normalized block sequence of any prescribed normalized disjoint sequence
such that every $T\in\mathcal A$ has a compact restriction to its span.
\end{corollary}
\begin{proof}
Apply Theorem~\ref{thm:countable} to a countable dense subset. Restriction
is norm-continuous and compact operators are closed in operator norm.
\end{proof}

\begin{proposition}[The separable-family conclusion is sharp]\label{prop:separable}
There is a norm-bounded, norm-separable family of compact operators on
$\ell_2$ for which no closed infinite-dimensional subspace makes all
restrictions nuclear.
\end{proposition}
\begin{proof}
Let $H=\ell_2$ and let $\mathcal A$ be the closed unit ball of
$\mathcal K(H)$, the space of compact operators. Finite matrices with
rational real and imaginary parts show that $\mathcal K(H)$ is separable
in operator norm; therefore so is $\mathcal A$. Every member is DSS.

\Needspace{7\baselineskip}
For any closed infinite-dimensional $F\subset H$, choose an orthonormal
basis $(w_n)$ of $F$ and define
\[
 S_Fx=\sum_{n\geq1}\frac{\langle x,w_n\rangle}{n}w_n.
\]
We take the inner product linear in its first variable. This operator is
zero on $F^\perp$, compact, and of norm one, so $S_F\in\mathcal A$.
If $S_F|_F=\sum_j f_j\otimes y_j$ were a nuclear representation, then
for every $m$ the Cauchy--Schwarz and Bessel inequalities would give
\begin{align*}
 \sum_{n=1}^m\frac1n
 &=\left|\sum_j\sum_{n=1}^m
      f_j(w_n)\langle y_j,w_n\rangle\right|\\
 &\leq\sum_j
   \left(\sum_{n=1}^m|f_j(w_n)|^2\right)^{1/2}
   \left(\sum_{n=1}^m|\langle y_j,w_n\rangle|^2\right)^{1/2}
 \leq\sum_j\norm{f_j}\norm{y_j}<\infty.
\end{align*}
The bound is independent of $m$, contradicting divergence of the harmonic
series. Thus $S_F|_F:F\to H$ is not nuclear.
\end{proof}

\begin{remark}
This is an obstruction to the conclusion itself, not just to the density
proof. The family in Proposition~\ref{prop:separable} is uncountable;
Theorem~\ref{thm:countable} still applies to any chosen countable subset.
For a norm-compact family the stronger uniform conclusion below holds.
\end{remark}

\section{Uniform nuclear tails at prescribed rates}\label{sec:uniform}

\begin{theorem}\label{thm:uniform}
Let $\mathcal K\subset\DSS(E,Y)$ be compact in operator norm. Fix a
normalized disjoint sequence $(x_n)$ and any positive null sequence
$(r_m)_{m\geq0}$. There is a normalized successive block sequence
$(v_n)$ of $(x_n)$ such that, for $F=[v_n]$ and its basis projections $P_m$,
\begin{equation}\label{eq:uniformnuclear}
 \sup_{T\in\mathcal K}\nucl\bigl(T|_F(I_F-P_m)\bigr)\leq r_m
 \quad(m\geq0).
\end{equation}
In particular, for the approximation numbers
$a_k(S)=\inf_{\operatorname{rank}R<k}\norm{S-R}$,
\[
 \sup_{T\in\mathcal K}a_{m+1}(T|_F)\leq r_m.
\]
The set $\{T|_F:T\in\mathcal K\}$ is compact in nuclear norm.
\end{theorem}
\begin{proof}
Set $\eta_n=2^{-n}\min_{0\leq k<n}r_k$. Then $\eta_n>0$ and
$\sum_{n>m}\eta_n\leq r_m$ for every $m\geq0$.
At stage $n$ choose a finite $\eta_n/2$-net for $\mathcal K$ with centers
in $\mathcal K$. Apply Lemma~\ref{lem:finite} to those centers on the
unused tail, obtaining a unit block $v_n$ on which every center has
image norm less than $\eta_n/2$. Therefore
\[
 \sup_{T\in\mathcal K}\norm{Tv_n}\leq\eta_n.
\]
The rank-one representations of Lemma~\ref{lem:nuclear} now give
\eqref{eq:uniformnuclear}. The operator $T|_FP_m$ has rank at most $m$,
which proves the approximation-number bound.

For fixed $m$, the map $T\mapsto T|_FP_m$ is continuous from operator norm
to nuclear norm, since
$\nucl((T-S)|_FP_m)\leq m\norm{T-S}$.
Its image of $\mathcal K$ is nuclear-norm compact. The uniform tail
estimate makes these images uniformly approximate the restricted family
in nuclear norm, proving its total boundedness. The restriction map itself
is continuous into nuclear norm by the same head-and-tail estimate;
its image of $\mathcal K$ is thus compact.
\end{proof}

\subsection*{A compact exhaustion of the operator family}

The countable and compact statements admit a joint extension. Here
\emph{norm-$\sigma$-compact} means a countable union of sets compact
in operator norm; replacing these sets by their finite unions makes
the exhaustion increasing. All statements below hold over either
scalar field and do not use Theorem~\ref{thm:main}.

\begin{theorem}[Simultaneous restrictions for a compact exhaustion]
\label{thm:sigma}
Let $\mathcal A\subset\DSS(E,Y)$ have an increasing compact exhaustion
$\mathcal A=\bigcup_{j\geq1}\mathcal K_j$ in operator norm. Fix a
normalized disjoint sequence $(x_n)$ and a positive null sequence
$(r_m)_{m\geq0}$. There are normalized successive blocks $(v_n)$ of
$(x_n)$, with $F=[v_n]$, such that $T|_F$ is nuclear for every
$T\in\mathcal A$ and
\begin{equation}\label{eq:sigmatail}
 \sup_{T\in\mathcal K_j}
 \nucl\bigl(T|_F(I_F-P_m)\bigr)\leq r_m
 \qquad(j\geq1,\ m\geq j-1).
\end{equation}
Moreover, there are nuclear operators $N_{T,m}:E\to Y$, for every
$T\in\mathcal A$ and $m\geq0$, such that, with $F_m=[v_n:n>m]$,
\begin{align}
 F_m&\subseteq\bigcap_{T\in\mathcal A}\ker(T-N_{T,m}),
                                                        \label{eq:sigmakernel}\\
 \sup_{T\in\mathcal K_j}\nucl(N_{T,m})&\leq r_m
 \qquad(m\geq j-1).                                    \label{eq:sigmacorrection}
\end{align}
For each fixed $j$ and $m$, both maps
$T\mapsto T|_F$ and $T\mapsto N_{T,m}$ are continuous from operator
norm on $\mathcal K_j$ to nuclear norm. Their images of $\mathcal K_j$
are therefore compact in nuclear norm.
\end{theorem}
\begin{proof}
Put $\eta_n=2^{-n}\min_{0\leq k<n}r_k$. At stage $n$, take a finite
$\eta_n/2$-net for $\mathcal K_n$ with centers in $\mathcal K_n$.
Lemma~\ref{lem:finite}, applied to those centers on an unused tail,
gives a normalized finite block $v_n$ with
\[
 \sup_{T\in\mathcal K_n}\norm{Tv_n}\leq\eta_n.
\]
If $n\geq j$, then $\mathcal K_j\subseteq\mathcal K_n$, so the
same estimate holds uniformly on $\mathcal K_j$. For $m\geq j-1$,
every $n>m$ satisfies $n\geq j$, and consequently
\begin{equation}\label{eq:sigmasum}
 \sup_{T\in\mathcal K_j}\sum_{n>m}\norm{Tv_n}
 \leq\sum_{n>m}\eta_n
 \leq r_m\sum_{n>m}2^{-n}\leq r_m.
\end{equation}
For each $T\in\mathcal K_j$ the first $j-1$ image norms form a
finite sum. Thus all its image norms are summable, proving nuclearity
of $T|_F$. Lemma~\ref{lem:nuclear} and \eqref{eq:sigmasum} give
\eqref{eq:sigmatail}.

Extend each coefficient functional $v_n^*$ from $F$ to
$\widetilde f_n\in E^*$ with norm one, making these Hahn--Banach
choices once for the entire family. Define
\[
 N_{T,m}x=\sum_{n>m}\widetilde f_n(x)Tv_n.
\]
The preceding summability makes this a nuclear operator for every
$T\in\mathcal A$. It agrees with $T$ on every $v_n$ with $n>m$,
giving \eqref{eq:sigmakernel}; \eqref{eq:sigmasum} gives
\eqref{eq:sigmacorrection}.

Finally, fix $j,m$ and let $\ell\geq\max(m,j-1)$. For either map
$\Phi(T)=T|_F$ or $\Phi(T)=N_{T,m}$, its rank-one representation
and \eqref{eq:sigmasum} give, for $T,S\in\mathcal K_j$,
\[
 \nucl\bigl(\Phi(T)-\Phi(S)\bigr)
 \leq\ell\norm{T-S}+2r_\ell.
\]
First choose $\ell$ to control the tail and then $\norm{T-S}$ to
control the finite head. This proves continuity and hence compactness
of both images.
\end{proof}

\begin{corollary}[Continuous parameter families]\label{cor:parameters}
Let $d\geq1$ and let $t\mapsto T_t$ be an operator-norm-continuous
map from $\mathbb R^d$ to $\DSS(E,Y)$. From any normalized disjoint
sequence one can select a single normalized block sequence $(v_n)$
such that all restrictions $T_t|_{[v_n]}$ are nuclear and the map
$t\mapsto T_t|_{[v_n]}$ is continuous in nuclear norm. Given any
positive null sequence $(r_m)$, the choice can also ensure the bound
\eqref{eq:sigmatail} uniformly on every compact set of parameters
for all sufficiently large $m$.
\end{corollary}
\begin{proof}
Apply Theorem~\ref{thm:sigma} to
$\mathcal K_j=\{T_t:t\in[-j,j]^d\}$. Every compact set of parameters
lies in one such cube, and every parameter has a neighbourhood inside
a sufficiently large cube. The uniform estimate and continuity on
each $\mathcal K_j$ therefore give the conclusions.
\end{proof}

\begin{remark}
The bounds in Theorem~\ref{thm:sigma} are uniform on each prescribed
compact piece after a finite initial segment. No uniform nuclear bound
over the whole union, and no global operator-norm-to-nuclear-norm
continuity on that union, is asserted. A bounded countable family can
already fail to have uniformly vanishing restriction tails; see
Proposition~\ref{prop:distinctions}. Conversely, norm separability alone
does not guarantee simultaneous nuclearity, by
Proposition~\ref{prop:separable}. Norm-$\sigma$-compactness is a sufficient
condition, not a claimed characterization of all admissible families.
The proof is a diagonal refinement of the preceding classical
small-restriction method; independent novelty has not been established.
\end{remark}

\section{Common Hilbert factors and small nuclear corrections}\label{sec:kernel}

The following refinements use the same block selection with faster image
decay. They apply to DSS operators over either scalar field, independently
of the candidate LSS results.

\begin{theorem}[Common factors and a common kernel]\label{thm:kernel}
Let $\mathcal K\subset\DSS(E,Y)$ be compact in operator norm, let
$(x_n)$ be normalized and disjoint, and let $\epsilon>0$. There exist
a normalized successive block sequence $(v_n)$ of $(x_n)$, $F=[v_n]$,
and nuclear operators $J:F\to\ell_2$ and $\widetilde J:E\to\ell_2$
such that $\widetilde J|_F=J$, $\nucl(\widetilde J)\leq1$, and $J$
is injective with dense range. Let $D:\ell_2\to\ell_2$ be the fixed
diagonal operator $De_n=2^{-n}e_n$. Then:
\begin{enumerate}
\item For every integer $k\geq0$ and $T\in\mathcal K$ there is a
nuclear $B_{T,k}:\ell_2\to Y$ with
\[
 T|_F=B_{T,k}D^kJ,\qquad
 \sup_{T\in\mathcal K}\nucl(B_{T,k})
 \leq\epsilon\sum_{n\geq1}2^{-n^2+(k+1)n}<\infty.
\]
For fixed $k$, $T\mapsto B_{T,k}$ is continuous from operator norm on
$\mathcal K$ to nuclear norm, and its image is compact.
\item There are nuclear $N_T:E\to Y$, $T\in\mathcal K$, such that
\[
 \sup_{T\in\mathcal K}\nucl(N_T)<\epsilon,
 \qquad F\subseteq\bigcap_{T\in\mathcal K}\ker(T-N_T).
\]
The map $T\mapsto N_T$ is continuous into nuclear norm and its image
is compact. One may take $N_T=B_{T,0}\widetilde J$.
\end{enumerate}
\end{theorem}
\begin{proof}
Use the finite-net selection in Theorem~\ref{thm:uniform}, now with
$\eta_n=\epsilon2^{-n^2}$, to choose $(v_n)$ satisfying
\begin{equation}\label{eq:fastimages}
 \sup_{T\in\mathcal K}\norm{Tv_n}\leq\epsilon2^{-n^2}.
\end{equation}
Extend $v_n^*\in F^*$ to $\widetilde f_n\in E^*$ with norm one by
Hahn--Banach, making these choices once for the whole family. Define
\[
 \widetilde Jx=\sum_{n\geq1}2^{-n}\widetilde f_n(x)e_n,
 \qquad J=\widetilde J|_F.
\]
This representation has nuclear sum $\sum_n2^{-n}=1$.
The formula $Jv_n=2^{-n}e_n$ proves injectivity on $F$ by uniqueness
of basis coefficients and proves density of its range.

For $b=(b_n)\in\ell_2$ define
\[
 B_{T,k}b=\sum_{n\geq1}b_n2^{(k+1)n}Tv_n.
\]
Equation~\eqref{eq:fastimages} makes this a nuclear representation with
the stated bound, since a quadratic negative exponent dominates any
fixed linear exponent. The identity
$B_{T,k}D^kJv_n=Tv_n$ holds on each basis vector and hence on $F$.
All powers use the same $F,J,D$; no uniform bound in $k$ is asserted.

Finally set
\begin{equation}\label{eq:correction}
 N_Tx=\sum_{n\geq1}\widetilde f_n(x)Tv_n.
\end{equation}
This is $B_{T,0}\widetilde J$, and the displayed representation gives
\[
 \nucl(N_T)\leq\epsilon\sum_{n\geq1}2^{-n^2}<\epsilon.
\]
It agrees with $T$ on each $v_n$ and therefore on $F$.

For continuity, put $a_n=2^{(k+1)n}$ for $B_{T,k}$ and $a_n=1$ for
$N_T$. The relevant representations give, for $T,S\in\mathcal K$,
\[
 \nucl(\Phi(T)-\Phi(S))
 \leq\left(\sum_{n\leq m}a_n\right)\norm{T-S}
       +2\epsilon\sum_{n>m}a_n2^{-n^2},
\]
where $\Phi(T)=B_{T,k}$ or $N_T$, respectively. First choose $m$ to
control the tail, then $\norm{T-S}$ to control the finite head.
Continuity and compactness of $\mathcal K$ prove both compact-image
assertions.
\end{proof}

\begin{remark}
Each $T$ has its own correction $N_T$; the block subspace, coordinate
extensions, and Hilbert factor are common. No bounded projection from $E$
onto $F$ has been constructed. In particular, $J$ is not claimed to be
bounded below and does not identify $F$ with a Hilbert space.
The assignment $T\mapsto N_T$ extends algebraically and linearly to
$\operatorname{span}\mathcal K$ via \eqref{eq:correction}, but no
operator-norm bounded extension to its norm closure is asserted.
\end{remark}

\begin{corollary}[A nuclear-perturbation characterization of DSS]\label{cor:characterization}
For a bounded $T:E\to Y$, the following are equivalent:
\begin{enumerate}
\item $T$ is DSS.
\item For every normalized disjoint sequence $(x_n)$ and every
$\epsilon>0$, there are normalized successive blocks $(v_n)$ of $(x_n)$
and a nuclear $N:E\to Y$ such that
\[
 \nucl(N)<\epsilon,\qquad [v_n]\subseteq\ker(T-N).
\]
\end{enumerate}
\end{corollary}
\begin{proof}
For the forward implication apply Theorem~\ref{thm:kernel} to
$\mathcal K=\{T\}$. Conversely, if $T$ is bounded below by $c>0$
on the span of a normalized disjoint sequence, choose $\epsilon<c$.
On the nonzero selected span one would have
\[
 c\norm v\leq\norm{Tv}=\norm{Nv}
 \leq\nucl(N)\norm v<c\norm v,
\]
a contradiction.
\end{proof}

\begin{proposition}[Finite-rank corrections do not suffice]\label{prop:finiterank}
There is a compact DSS operator $T:\ell_2\to\ell_2$ such that
$\ker(T-R)$ is finite-dimensional for every finite-rank $R$.
\end{proposition}
\begin{proof}
Take $Te_n=n^{-1}e_n$. This compact operator is injective. On
$\ker(T-R)$ it agrees with $R$, so it maps that kernel injectively
into $\operatorname{ran}R$. Hence
$\dim\ker(T-R)\leq\operatorname{rank}R$.
\end{proof}

\begin{corollary}[LSS consequences of the proposed main theorem]\label{cor:conditional}
If Theorem~\ref{thm:main} is valid, Theorems~\ref{thm:countable},
\ref{thm:uniform}, \ref{thm:sigma}, and \ref{thm:kernel}, and
Corollaries~\ref{cor:separable} and \ref{cor:parameters}, apply to real LSS families without
common-factor assumptions. The perturbation property of
Corollary~\ref{cor:characterization} then characterizes real LSS operators.
\end{corollary}
\begin{proof}
Use the asserted equality to replace LSS by DSS in those statements.
\end{proof}

\section{Scope, sharpness, and remaining questions}\label{sec:limits}

\begin{proposition}[Two necessary distinctions]\label{prop:distinctions}
In Theorem~\ref{thm:uniform}, one cannot generally require the selected
vectors to be positive.
Nor can one replace norm compactness in Theorem~\ref{thm:uniform} by
boundedness, even for a countable family of rank-one operators.
\end{proposition}
\begin{proof}
The functional $T:\ell_1\to\mathbb R$, $Tx=\sum_nx_n$, has rank one and
is therefore DSS. Yet $\norm{Tu}=1$ for every positive unit vector $u$.
In particular, the conclusion at $m=0$ with $r_0<1$ cannot hold on the
span of any normalized positive sequence.

For the second assertion, let $(h_j)$ be a countable dense subset of the
unit sphere of real $\ell_2$ and put $T_jx=\langle x,h_j\rangle$.
Each $T_j$ has rank one and norm one, but
$\sup_j|T_jx|=\norm x$ for every $x$.
On any infinite-dimensional block span $F$ and every $m$, taking a unit
vector in its tail gives
\[
 \sup_j\norm{T_j|_F(I_F-P_m)}=1.
\]
Uniform tails tending to zero are impossible. Individual nuclearity does
not imply uniform smallness of the whole family.
\end{proof}

The subspace $F$ in Theorem~\ref{thm:uniform} depends on the family, the
initial disjoint sequence, and the prescribed rates. Arbitrarily fast
decay on this selected subspace is not a universal quantitative comparison
of LSS and DSS for the original operator.

The DSS selection results are valid over the complex field as well.
This does not establish complex LSS$=$DSS. The real alternating-chain
proof uses positive and negative parts of signed coefficients. A lower
bound for a real restriction need not survive on its complex span:
the map $\mathbb C^2\to\mathbb C$, $(z_1,z_2)\mapsto z_1+iz_2$,
is isometric on $\mathbb R^2$ with its Euclidean norm but has a complex
kernel. This finite-dimensional example only blocks a naive transfer;
it is not an LSS--DSS counterexample.

A quantitative comparison formulated directly for the original operator
would require control of the auxiliary quotient cone constant.
Corollary~\ref{cor:threshold} does not provide that control.

\begin{remark}[Field and quantifiers]
The proof splits real signed coefficients into positive and negative
parts. No complex-field extension is asserted. The conclusion covers
arbitrary real source lattices and real Banach-space targets:
reflexivity and the unconditional basis are properties of the auxiliary
witness space, not hypotheses on $E$ or $X$.
The chain length is chosen after $c,\delta,\norm B$ and $\norm C$
have been fixed. The row limit is then taken at that fixed finite
length. There is no exchange of these two quantifiers.
\end{remark}

\begin{remark}[External inputs and points for review]
The external inputs to the main argument are the disjoint $c_0/\ell_1$
exclusion in \cite{FLT} and the unconditional-basis reflexivity criterion.
Theorem 2.5 of \cite{FLT} is used only after the main theorem, for the
sequence formulation and the interpretation of additivity. Independent
review should focus on Lemmas~\ref{lem:quotient} and \ref{lem:chains},
especially the common coefficients and normalization
\eqref{eq:normalization}. The finite estimate
Corollary~\ref{cor:threshold} can be checked separately.
\end{remark}

\subsection*{Use of generative AI tools}
OpenAI ChatGPT was used extensively for literature exploration,
development and internal checking of mathematical arguments, and
preparation and editing of this manuscript. Such assistance does not
constitute independent verification. Ma Yicen is the sole named author;
responsibility for the contents of a submitted version rests with the
author.

\end{document}